\documentclass[12pt]{article}

\usepackage{amsmath}
\usepackage{amsfonts}
\usepackage{amssymb}
\usepackage{amscd}

\newtheorem{Lemma}{Lemma}[section] \newtheorem{Theorem}[Lemma]{Theorem}
\newtheorem{Proposition}[Lemma]{Proposition}
 
\newtheorem{Definition}[Lemma]{Definition}

\def\vfi{\varphi}
 
\def\bC{{\mathbb C}} 
\def\bT{{\mathbb T}}
\def\H{{\cal H}} 
\def\A{{\cal A}}
\def\C{{\cal C}}
\def\S{{\cal S}} 
\def\U{{\cal U}}
\def\Q{{\cal Q}}
\def\F{{\cal F}}
\def\M{{\cal M}}
\def\Op{\mathfrak{Op}} 
\def\Rep{\mathfrak{Rep}} 

\def\h{\hbar}
\def\gA{\mathfrak A}
\def\gB{\mathfrak B}
\def\gC{\mathfrak C}
\def\gF{\mathfrak F}

\def\n{\parallel}
\def\bR{\mathbb R}
\def\bN{\mathbb N}

\def\bF{\mathbb F}
\def\pl{\parallel}
\def\pt{\partial}
\def\n{|\!|\!|}

\begin{document}

\title{Strict Deformation Quantization for\\ a Particle in a Magnetic Field}

\date{\today}

\author{Marius M\u antoiu and Radu Purice 
\footnote{Institute of Mathematics ``Simion Stoilow'' of
the Romanian Academy, P.O.  Box 1-764, Bucarest, RO-014700, Romania, Electronic 
mail: Marius.Mantoiu@imar.ro, Radu.Purice@imar.ro}}

\maketitle

\begin{abstract} 

Recently, we introduced a mathematical framework for the quantization of a particle in a variable magnetic field. It consists in a modified form of the Weyl pseudodifferential calculus and a $C^*$-algebraic setting, these two points of view being isomorphic in a suitable sense. In the present paper we leave Planck's constant vary, showing that one gets a strict deformation quantization in the sense of Rieffel. In the limit $\h\rightarrow 0$ one recovers a Poisson algebra induced by a symplectic form defined in terms of the magnetic field. 

\end{abstract}



\section*{Introduction} 

The present article treats the semiclassical limit of the mathematical formalism describing a quantum, non-relativistic particle without internal structure, placed in a variable magnetic field. The limit is considered in the precise sense of Rieffel's axioms (cf. \cite{Ri2}, \cite{Ri3}, \cite{La}), involving $C^*$-algebras. This setting is widely called {\it strict} quantization, to distinguish it from the version in terms of formal series (see \cite{BF} for example). It consists of several ingredients, which we outline here very briefly, refering to Section 1 for a detailed discussion: 

1. One needs first a natural family of classical observables. It is admitted that this should form a Poisson algebra $\gA$, which is roughly a real associative and commutative algebra endowed with a compatible Poisson bracket. This structure describes the classical physical system.

2. For non-null values of Planck's constant $\h$, one has to define $C^*$-algebras of quantum observables $\gC^\h$.

3. It must be shown that for $\h\rightarrow 0$ ``the quantum structure converges to the classical one''. This is described precisely by Rieffel's system of axioms or some of its versions.

If a certain extra technical condition is verified, allowing to define on classical observables a family of ``deformed products'' indexed by $\h$, one speaks of {\it strict deformation quantization}.

In our case, a particle without spin moving in the $N$-dimensional
configuration space $\bR^N$ and placed under the influence of an external
variable magnetic field, the natural Poisson algebra is well-known. The
observables are smooth functions defined on the phase space $\Xi:=\bR^{2N}$,
the associative product is defined pointwise and the Poisson bracket is
induced by the canonical symplectic form on $\Xi$, to which we add a magnetic
contribution (\cite{KO1}, \cite {MR}). This is described in Section 2. 

Quite surprisingly, the algebras of quantum observables for this system were defined and developed only recently. One reason could be that the canonical variables in the magnetic case (the components of the position and those of the magnetic momentum) satisfy complicated commutation relations, that have to be taken into account when defining more general observables as functions of these basic ones. The intensive use of {\it constant} magnetic fields and (or) special observables that are {\it quadratic} with respect to the momenta have also played a certain role. The setting which is correct (at least in our opinion) appeared in \cite{KO1}, \cite{MP2} and \cite{KO2} (a pseudodifferential point of view) and in \cite{MP1} and \cite{MPR1} ($C^*$-algebras). The right attitude can also be found in \cite{Lu}, but undeveloped and stated for a very particular case; it seems that it has been largely unnoticed. The critical point is {\it gauge invariance}: when several equivalent vector potentials corresponding to a given magnetic field are used in defining observables, the results should be connected by simple unitary equivalences. But to achieve this, one has to be very careful in defining the precise form of the observables as well as the composition laws to which they are submitted. We shall explain all these in Section 3. 

In Section 4 we state our Main Theorem. It asserts that under certain hypothesis, to a magnetic field and to an abelian algebra of ``configurational observables'' one can associate naturally a strict deformation quantization.

Sections 5, 6 and 7 are devoted to the proof of the Main Result. The three non-trivial axioms are verified separately. The Rieffel condition and the von Neumann condition follow from the 
results of \cite{PR2} and \cite{Ni} (Sections 5 and 6). Nevertheless, we also give a direct elementary proof for the von Neumann condition similar with that given for the Dirac condition (Sections 6 and 7).

This article is addressed also to people that have not deformation
quantization as their main skill. The system we treat has a certain physical
interest (this is not always the case in this field). Thus we decided to avoid
technical complications and to leave more general situations to subsequent
works. In particular, we hope to be able to say something on strict
deformation quantization by twisted groupoids (see \cite{La}, \cite{LR}, \cite{LMN} \cite{NWX}, and
references therein for the untwisted case), which should include the present
work as a particular instance. A pure state quantization would also be an
interesting topic.

Both in classical and in quantum theory one works with ``real''
observables. For any space $\mathcal E$ of complex functions we denote by
$\mathcal E_\mathbb R$ the subspace of $\mathbb R$-valued elements in
$\mathcal E$. For instance, $C^\infty(P)_\mathbb R$ will be the family of real
$C^\infty$ functions on the smooth manifold $P$. If $\mathfrak C$ is a
$C^*$-algebra, we set $\mathfrak C_\mathbb R$ for the set of self-adjoint
elements of $\mathfrak C$.

Some other notations: If $Y$ is a locally compact group we denote by $C(Y)$
 the $*$-algebra of all continuous complex functions on $Y$.  
 $BC(Y)$, $BC_u(Y)$, $C_0(Y)$ mean respectively ``bounded and continuous'',
 ``bounded and uniformly continuous'' and ``continuous and small at
 infinity''.  If $\H$ is a Hilbert space, $K(\H)$ will be the set of all
 compact operators in $\H$, forming an ideal in the $C^*$-algebra $B(\H)$ of
 all the linear, bounded operators on $\H$. The unitary operators form the
 group $\U(\H)$.

{\bf Acknowledgements:} A large part of this work has been completed while the authors visited the University of Geneva and we express our gratitude to Prof. Werner Amrein for his kind hospitality and the stimulating discusions. We also acknowledge the partial support from the EURROMMAT Programme (contract no. ICA1-CT-2000-70022) and from the CERES Programme (contract no. 3-28/2003).
We are greatful to Fr\'ed\'eric Cadet for a useful discussion and to Serge
Richard for a critical reading of the manuscript.

\section{The axioms}

We describe here Rieffel's framework for strict quantization. There are several versions of his axioms; we choose to work with the system of axioms which appears in \cite{La}, to which we also refer for many other details. The starting point is a ``classical algebra of observables'' described by a Poisson algebra.

\begin{Definition}\label{Landsman}
{\rm A Poisson algebra} is a triple $\left(\gA,\circ,\{\cdot,\cdot\}\right)$,
where $\gA$ is a real vector space, $\ \circ,\ \{\cdot,\cdot\}$ are bilinear
maps $:\gA\times\gA\rightarrow\gA$ such that $\circ$ is associative and commutative, $\{\cdot,\cdot\}$ is antisymmetric and for each $\ \vfi\in\gA$, $\{\vfi,\cdot\}$ is a derivation both with respect to $\circ$ and to $\{\cdot,\cdot\}$.
Thus, aside bilinearity, the two maps satisfy for all $\vfi,\psi,\rho\in \gA$:

(i) $\psi\circ\vfi=\vfi\circ\psi,\ \ \ $ $(\psi\circ\vfi)\circ\rho=\psi\circ(\vfi\circ\rho)$,

(ii) $\{\psi,\vfi\}=-\{\vfi,\psi\}$,

(iii) $\{\vfi,\psi\circ\rho\}=\psi\circ\{\vfi,\rho\}+\{\vfi,\psi\}\circ\rho\ \ $ {\rm (Leibnitz rule)},

(iv) $\{\vfi,\{\psi,\rho\}\}=\{\{\vfi,\psi\},\rho\}+\{\psi,\{\vfi,\rho\}\}\ $ {\rm (Jacobi's identity)}.
\end{Definition}

The elements of $\gA$ are interpreted as observables of a classical description of a physical system. For each $\vfi\in\gA$ and each value $\h\ne 0$ of Planck's constant, one would like to have an object $\Q^\h(\vfi)$ representing the same observable in a quantum description of the system. One also hopes that the algebraic structure of the quantum observables should converge to the classical picture described by the Poisson algebra, in some suitable norm $\pl\cdot\pl_\h$ depending continuously of $\h$. This might be seen as a precise mathematical form of Bohr's correspondence principle.

A systematic justification of the next definitions may be found in \cite{La}. Note that usually in $\gA$ many classical observables are ``unbounded''; the use of norms forces us to apply quantization only to certain subfamilies $\gA_0$ of $\gA$. 

\begin{Definition}\label{strict}
Let $\gA_0$ be a Poisson algebra which is densely contained in the 
self-adjoint part $\gC^0_\bR$ of an abelian $C^*$-algebra $\gC^0$.
{\rm A strict quantization of the Poisson algebra} $\left(\gA_0,\circ,\{\cdot,\cdot\}\right)$ is a family of maps $\left(\Q^\h:\gA_0\rightarrow\gC^\h_{\bR}\right)_{\h\in I}$, where

(i) $I$ is a subset of the real axis, for which the origin is an accumulation point contained in $I$,

(ii) $\gC^\h$ is a $C^*$-algebra, with product and norm denoted respectively by $\sharp^\h$ and $\parallel\cdot\parallel_\h$. For $\vfi^\h,\psi^\h\in\gC^\h_{\bR}$ (the self-adjoint part of $\gC^\h$) we set $\vfi^\h\star^\h\psi^\h:=\frac{1}{2}\left(\vfi^\h\sharp^\h\psi^\h+\psi^\h\sharp^\h\vfi^\h\right)$ {\rm (a Jordan product)} and $\{\vfi^\h,\psi^\h\}_\h:=\frac{1}{i\h}\left(\vfi^\h\sharp^\h\psi^\h-\psi^\h\sharp^\h\vfi^\h\right)$.

(iii) $\Q^\h:\gA_0\rightarrow\gC^\h_{\bR}$ is $\bR$-linear for each $\h$ and $\Q^0$ is just the inclusion map,

and the following axioms are fulfilled:

(a) RIEFFEL'S CONDITION: For $\ \vfi\in \gA_0$, the map $I\ni\h\rightarrow\parallel\Q^\h(\vfi)\parallel_\h\ \in\bR_+$ is continuous.

(b) VON NEUMANN CONDITION: For $\ \vfi,\psi\in \gA_0$, $\ \parallel\Q^\h(\vfi)\star^\h\Q^\h(\psi)-\Q^\h(\vfi\circ\psi)\parallel_\h\rightarrow 0$ when $\h\rightarrow 0$.

(c) DIRAC'S CONDITION: For $\ \vfi,\psi\in \gA_0$, $\ \parallel\{\Q^\h(\vfi),\Q^\h(\psi)\}_{\h}-\Q^\h\left(\{\vfi,\psi\}\right)\parallel_\h\rightarrow 0$ when $\h\rightarrow 0$.

(e) COMPLETENESS: $\ \Q^\h(\gA_0)$ is dense in $\gC^\h_{\bR}$ for all $\h\in I$.
\end{Definition}

The word ``strict'' was coined by Rieffel in order to distinguish his framework from the (deformation) quantization defined in terms of formal series.
Usually Poisson algebras are function spaces:

\begin{Definition}\label{vPoisson}
We call {\rm Poisson manifold} a smooth manifold $M$ so that on $C^\infty(M)_\bR$ a bracket $\{\cdot,\cdot\}$ is given such that, denoting by $\circ$ the pointwise multiplication, the triple $\left(C^\infty(M)_\bR,\circ,\{\cdot,\cdot\}\right)$ is a Poisson algebra.
\end{Definition}
 
When $M$ is not compact, $C^\infty(M)_\bR$ is a very large, unnormed space. In quantization one deals with suitable families of smooth bounded observables:

\begin{Definition}\label{strictQvP}
{\rm A strict quantization of the Poisson manifold} $M$ means the choice of a
Poisson subalgebra $\gA_0$ of $C^\infty(M)_\bR$ composed of bounded functions and a strict quantization of this Poisson subalgebra.
\end{Definition}

One should be aware that the linear maps $\Q^\h$ tend to behave as morphisms only in the asymptotic limit $\h\rightarrow 0$. But under favorable circumstances (fulfilled rather often, but by no means always) they may serve to define modified products on $\gC^0$. In this case, one really is allowed to think in terms of ``deformed products''.

\begin{Definition}\label{strictdefQ}
A strict quantization $\left(\Q^\h:\gA_0\rightarrow\gC^\h_{\bR}\right)_{\h\in I}$ is called {\rm a strict deformation quantization} if for each $\h$, $\Q^\h\left(\gA_0\right)$ is a subalgebra of $\gC^\h_\bR$ and $\Q^\h$ is injective.
\end{Definition}

In such a case, for any $\h$, one defines $\ \sharp^\h:\gA_0\times\gA_0\rightarrow\gA_0$ such that $\Q^\h(\vfi\sharp^\h\psi)=\Q^\h(\vfi)\sharp^\h\Q^\h(\psi)$ for all $\vfi,\psi\in\gA_0$. The notational ambiguity is deliberate.

\textbf{Remark.} At the suggestion of the referee, to whom we thank for his interesting observations, we shall briefly comment upon an alternative definition for \textit{strict deformation quantization of a Poisson algebra} $\gA_0$, in the spirit of Rieffel's approach \cite{Ri1}, \cite{Ri2}, \cite{Ri3}.  The maps $\Q^\h$ being injective (see Definition \ref{strictdefQ}), we may identify all the algebras $\Q^\h[\gA_0]$ and consider the different $C^*$-algebras $\gC^\h$ as completions for different $C^*$-norms $\|.\|_\h$ of the same $*$-algebra $\mathbb{C}\otimes\gA_0$. We denote by $\overline{\mathbb{C}\otimes\gA_0}$ the completion taken with respect to the $C^*$-norm $\n .\n :=\underset{\h\in I}{\sup}\|.\|_\h$. Then we may define a strict deformation quantization of the Poisson algebra $\gA_0$, as
a family $\{\gC^\h\}_{\h\in I}$ of $C^*$-algebras (with products $\sharp^\h$ and $C^*$-norms $\|.\|_\h$) such that:
\begin{itemize}
\item $\gC^0$ is abelian;
\item $\mathbb{C}\otimes\gA_0$ is dense in $\gC^\h$ for any $\h\in I$;
\item the triple $(I,\{\gC^\h\}_{\h\in I},\Gamma)$, with $\Gamma:=C(I;\overline{\mathbb{C}\otimes\gA_0})$, defines a 
\textit{a continuous field of $C^*$-algebras}, cf. \cite{Di};   
\item THE DIRAC CONDITION: For any $\varphi$ and $\psi$ in $\gA_0$ we have
$$
\underset{\h\rightarrow 0}{\lim}\left\|\frac{1}{i\h}\left(\varphi\sharp^\h\psi-\psi\sharp^\h\varphi\right)-
\{\varphi,\psi\}\right\|_\h=0.
$$
\end{itemize} 
In fact our proof in this paper may be seen to give such a structure.

\section{The magnetic Poisson algebra}

For the convenience of the reader, we start by recalling briefly the way a symplectic manifold acquires a canonical Poisson structure.
For a differentiable manifold $M$ we denote by $C^\infty(M)$ the vector space of smooth real functions on $M$, by $\mathcal X(M)$ the $C^\infty(M)$-module of vector fields on $M$ and by $\Omega^k(M)$ the $C^\infty(M)$-module of k-forms on $M$ (i.e. $C^\infty$ sections of the fibre bundle of antisymetric k-linear forms $\Lambda^k_mM$ on $\mathbb{T}_mM$, the tangent space of $M$ at $m\in M$). One has $\Omega^0(M)=C^\infty(M)$. We denote by $d:\Omega^k(M)\rightarrow\Omega^{k+1}(M)$ the exterior differential. A symplectic form on $M$ is just a closed nondegenerate 2-form $\Sigma\in\Omega^2(M)$.

It follows easily from the axioms that the bracket $\{\cdot,\cdot\}$ of any Poisson manifold $M$ is given by {\it a Poisson bivector}. This means that one has $\{f,g\}=w(df,dg)$ for all $f,g\in C^\infty(M)$, where $w:\Omega^1(M)\times\Omega^1(M)\rightarrow C^\infty(M)$ is bilinear, antisymmetric and satisfies an extra condition connected to the Jacobi identity (see \cite{La} or \cite{Va} for details). The symplectic form will lead to such a Poisson bivector in a specific way. Being a nondegenerate bilinear form on each tangent space, $\Sigma$ defines a $C^\infty(M)$-linear isomorphism $\beta:\Omega^1(M)\rightarrow\mathcal X(M)$, $\theta(v):=\Sigma_m(\beta_m(\theta),v),\;\forall(\theta,v)\in\mathbb{T}^*_mM\times\mathbb{T}_mM$. Then one defines 
\begin{equation}\label{Vaisman}
\{f,g\}_\Sigma :=\Sigma\left[\beta(df),\beta(dg)\right],
\end{equation}
so in this case the Poisson bivector is given by $w_\Sigma=\Sigma\circ(\beta\times\beta)$. One checks easily that, in this way, $M$ becomes a Poisson manifold and we denote by $\mathcal{P}_0(M)$ the algebra $C^\infty(M)_{\mathbb{R}}$ endowed with the pointwise multiplication and the above canonical Poisson bracket.

We come back to our specific situation. The configuration space of our particle without internal structure is the space $X:=\bR^N$, with elements $q,x,y,z$. The subsequent presence of a magnetic field demands $N\ge 2$. We denote by $X^\star$ the dual of the vector space $X$, with elements $p,k,l$ and by $(x,p)\mapsto x\cdot p$ the duality between $X$ and $X^\star$. 

The phase-space of the system is the cotangent bundle $\mathbb{T}^*X$ of $X$, often denoted by $\Xi$ and identified with the direct sum $X\times X^\star$ (by identifying all the fibres with $X^\star$, using the action through translations). Typical vectors in $\Xi$ are $\xi=(q,p)$ or $\eta=(x,k)$. All the tangent spaces $\mathbb{T}_\xi(\Xi)$ will be identified with $\Xi$ and all the cotangent spaces $\mathbb{T}^*_\xi(\Xi)$ will be viewed as $\Xi^\star$ and, furthermore, as $X^\star\times X$. On $\Xi$ we have the canonical (constant) symplectic form defined by
\begin{equation*}
\sigma:\Xi\times \Xi\rightarrow \bR,\ \ \ \sigma[(x,k),(y,l)]:=y\cdot k-x\cdot l.
\end{equation*} 
This structure is adequate for the description of the particle whith no magnetic field. When a magnetic field is present, this can be taken into account by a change in the symplectic structure, cf. \cite{MR}. 

We thus consider a special class of flat symplectic manifolds, representing 'perturbations' of the above symplectic space and associated to a general (regular) magnetic field on $X$. In fact such \textit{a magnetic field is described by a closed 2-form} $B\in\Omega^2(X)$. Starting with the canonical projection $\pi:\Xi\cong X\times X^\star\rightarrow X$, we define canonically an injection $\tilde{\pi}_2:\Omega^2(X)\rightarrow\Omega^2(\Xi)$. Thus we get a new symplectic form $\sigma_B$ on $\Xi$ as the sum $\sigma_B:=\sigma+\tilde{\pi}_2B$, i.e.
$$
\left(\sigma_B\right)_{(q,p)}[(x,k),(y,l)]:=\sigma[(x,k),(y,l)]+(\tilde{\pi}_2B)_{(q,p)}[(x,k),(y,l)]=
$$
$$
=y\cdot k-x\cdot l + B_{q}(x,y).
$$
Being the sum of two closed forms, this 2-form is closed. It is also nondegenerate, thus it is a symplectic form on $\Xi$. Then (\ref{Vaisman}) gives
\begin{equation*}
\{f,g\}_B\equiv\{f,g\}_{\sigma_B}=\sigma_B\left[\beta(df),\beta(dg)\right],
\end{equation*}
so we badly need an explicit formula for $\beta$. Let us denote by
$\left<\cdot,\cdot\right>$ the duality between $\Xi$ and $\Xi^\star$. The
inverse $\beta^{-1}_\xi:\Xi\rightarrow\Xi^\star$ is defined, for $\xi=(q,p),
\eta=(x,k), \zeta=(y,l)\in\Xi$, by 
\begin{equation*}
(\sigma_B)_{(q,p)}[(x,k),(y,l)]\equiv y\cdot k-x\cdot l+x\cdot\overline{B}_q y=\left<(x,k),\beta^{-1}_{(q,p)}(y,l)\right>, 
\end{equation*}
where $\overline{B}_q:X\rightarrow X^\star$ is the linear, antisymmetric operator defined by $B_q(x,y)=x\cdot\overline{B}_q y$, $\forall x,y\in X$. It follows easily that $\beta^{-1}$ can be put in matrix-form
\begin{equation*}
\beta^{-1}_\xi=\left(
\begin{array}{cc}
\overline{B}_q&-1_{X^*}\\
1_{X}&0
\end{array}
\right)
\;:\;X\times X^\star\rightarrow X^\star\times X, 
\end{equation*}
which leads to the next matrix-form of $\beta_\xi$:
\begin{equation*}
\beta_\xi=\left(
\begin{array}{cc}
0&1_{X}\\
-1_{X^*}&\overline{B}_q
\end{array}
\right)
\ :\;X^\star\times X\rightarrow X\times X^\star.
\end{equation*}
Thus, writting $d_\xi h=\left(d^X_\xi h,d^{X^\star}_\xi h\right)$, one gets 
\begin{equation*}
\{f,g\}_B(\xi)=d^X_\xi f\cdot d^{X^\star}_\xi g-d^X_\xi g\cdot d^{X^\star}_\xi f+B_q\left(d^X_\xi f,d^X_\xi g\right).
\end{equation*}
Using coordinates, one has $d^X_\xi
h=\sum_{j=1}^N\left(\partial_{p_j}h\right)(\xi)\ dp_j$ and $d^{X^\star}_\xi
h=\sum_{j=1}^N\left(\partial_{q_j}h\right)(\xi)\ dq_j$ (recall that
$dp_j\in\left(X^\star\right)^\star\equiv X$). We get finally 
\begin{equation}\label{poisson}
\{f,g\}_B=\sum_{j=1}^N\left(\partial_{p_j}f\ \partial_{q_j}g-\partial_{q_j}f\ \partial_{p_j}g\right)+\sum_{j,k=1}^N B_{jk}(\cdot)\ \partial_{p_j}f\ \partial_{p_k}g.
\end{equation}

We shall denote by $\mathcal{P}_B(\Xi)$ the Poisson algebra
$C^\infty(\Xi)_{\mathbb{R}}$ endowed with the pointwise multiplication and the Poisson bracket $\{\cdot,\cdot\}_B$ given in (\ref{poisson}). The different descriptions of the quantum observable algebras in the next Section asks also for a partial Fourier transformed version of this Poisson algebra. This will be explained in Section 4 under favorable circumstances.

\section{The quantum magnetic observables; \\
 the $C^*$-algebras for $\hbar\ne 0$}

We are placed in the framework of the previous Section, but for most of the constructions the smoothness assumption on $B$ will be useless. We assume for the moment only that $B$ is continuous.

We present first a pseudodifferential approach to the magnetic quantum system, following \cite{MP2}, \cite{KO1} and \cite{KO2}. No $C^*$-algebras are in sight for the moment. We choose some vector potential $A$ corresponding to the magnetic field $B$ ($dA=B$). It also can be chosen continuous; think of the transversal gauge for example. The vector potential is used to define a representation of some explicitely gauge invariant structure. Only this one will be used in the process of quantization. 

Let us fix some value $\h\ne 0$ for the Planck constant. We would like to justify the construction of a correspondence $f\mapsto \Op
^\h_A(f)$ between (suitable) complex functions defined on the phase space $\Xi$ and operators. To the function $(q,p)\mapsto q_j$ one wants to assign the operator $Q_j$ of multiplication with $q_j$ (i.e.: $(Q_ju)(q):=q_ju(q)$) and to $(q,p)\mapsto p_j$ we associate the first-order differential operator $\Pi^\h_{A,j}:=\h P_j-A_j(Q)=-i\h\partial_j-A_j$. The difficulty of defining a functional calculus $f\mapsto \Op^\h_A(f)\equiv f(Q,\Pi_A^\h)$ for these $2N$ operators comes from their high degree of non-commutativity:
\begin{equation*}
i[Q_j,Q_k]=0,\ \ \ i[\Pi^\h_{A,j},Q_k]=\h\delta_{j,k},\ \ \ i[\Pi^\h_{A,j},\Pi_{A,k}^\h]=\h B_{kj}(Q),\ \ \ j,k=1,\dots,N.
\end{equation*}
A convenient global form of these canonical commutation relations may be given in terms of {\it the magnetic Weyl system}. Recall the unitary group $\left(e^{iQ\cdot p}\right)_{p\in X^\star}$ of the position as well as {\it the magnetic translations} $\left(U_A^\h(q):=e^{iq\cdot\Pi_A^\h}\right)_{q\in X}$, given explicitely in the Hilbert space $\H:=L^2(X)$ by
\begin{equation}\label{uuu}
U_A^\h(q)=e^{-(i/\h)\Gamma_A([Q,Q+\h q])}e^{iq\cdot\h P},
\end{equation}
where $\Gamma_A([q',q'+\h q]):=\int_{[q',q'+\h q]}A$ is the circulation of the vector potential $A$ along the segment 
$$
[q',q'+\h q]:=\{q'+t\h q\mid t\in[0,1]\}.
$$
The family $\left(U_A^\h(q)\right)_{q\in X}$ satisfies
\begin{equation*}
U_A^\h(q)U_A^\h(q')=\omega^\h_B(Q;q,q')U_A^\h(q+q'),\ \ \ q,q'\in X,
\end{equation*}
where we set 
\begin{equation*}
\omega^\h_B(q_0;q,q'):=e^{-(i/\h)\Gamma_B(<q_0,q_0+\h q,q_0+\h q+\h q'>)}
\end{equation*} 
and  
\begin{equation*}
\Gamma_B(<q_0,x,y>):=\int_{<q_0,x,y>}B
\end{equation*}
is the flux of $B$ through the triangle $\ <q_0,x,y>\ $ defined by the points $q_0,x$ and $y$.

Now the magnetic Weyl system is the family $\left(W_A^\h(q,p)\right)_{(q,p)\in\Xi}$ of unitary operators in $\H$ given by
\begin{equation*}
W_A^\h(q,p):=e^{-i\sigma\left((q,p),(Q,\Pi_A^\h)\right)}=e^{-i(Q+(\h/2)q)\cdot p}e^{-(i/\h)\Gamma_A([Q,Q+\h q])}e^{iq\cdot\h P}
\end{equation*}
and it satisfies for all $\ (q,p),(q',p')\in\Xi$
\begin{equation*}
W_A^\h(q,p)W_A^\h(q',p')=e^{(i/2)\sigma\left((q,p),(q',p')\right)}\omega_B^\h(Q;q,q')W_A^\h(q+q',p+p').
\end{equation*}

To construct $\Op^\h_A(f)\equiv f(Q,\Pi_A^\h))$ one does not dispose of a
spectral theorem. Having the functional calculus with a $C_0$-group in mind
and having faith in the ability of the magnetic Weyl system to take into
account the way $(Q_1,\cdots,Q_N;\Pi^\h_{A,1},\cdots,\Pi^\h_{A,N})$ fail to
commute, one proposes 
\begin{equation*}
\Op^\h_A(f):=\int_\Xi d\xi\left(\gF_\Xi f\right)(\xi)W_A^\h(\xi),
\end{equation*}
where by $\gF_\Xi$ we denote the symplectic Fourier transform
\begin{equation*}
\left(\gF_\Xi f\right)(\xi):=\int_\Xi d\eta\ e^{-i\sigma(\xi,\eta)}f(\eta).
\end{equation*}
A suitable choice of the Haar measures on $X$, $X^\star$ and $\Xi$ leads to the exact form of the formulae above, with no numerical factors in front of the integrals.

Some simple replacements lead to the following expression for the action of these operators on vectors $u\in L^2(X)$:
\begin{equation}\label{Op}
\left[\Op^\h_A(f)u\right](x)=\h^{-N}\int_X\int_{X^\star}dy\ dk \ e^{(i/\h)(x-y)\cdot k}e^{-(i/\h)\Gamma_A([x,y])}f\left(\frac{x+y}{2},k\right)u(y).
\end{equation}
To have $\ \Op^\h_A(f)\Op^\h_A(g)=\Op^\h_A(f\circ^\h g)$ and $\
\Op^\h_A(f)^*=\Op^\h_A(f^{\circ^\h})$, one sets
$f^{\circ^\h}(q,p):=\overline{f(q,p)}$ (independent of $\h$ or $B$) and 
\begin{equation}\label{circ}
(f\circ^\h g)(\xi):=(2/\h)^{2N}\int_{\Xi}d\eta\int_{\Xi}d\zeta \ e^{-2(i/\h)\sigma(\xi-\eta,\xi-\zeta)}e^{-(i/\h)\Gamma_B(<q-y+x,x-q+y,y-x+q>)}f(\eta)g(\zeta).
\end{equation}
The composition law $\circ^\h\equiv\circ^\h_B$ depends only on the magnetic field and not on the choice of some vector potential.

Obviously, for $B=0$ and $A=0$ the above formulae reproduce the well-known formulae of the pseudodifferential calculus in Weyl form. In \cite{MP2} they are studied in detail, their gauge-invariance is underlined and a rigorous meaning of them and of some of their extensions are outlined. See also \cite{KO1}, \cite{KO2} for other developments and for nice geometrical interpretations. We shall come back to this {\it magnetic Weyl calculus} after an excursion into twisted crossed product algebras.

The input for a crossed product is a locally compact group $X$ acting on a $C^*$-algebra $\A$. One constructs a larger $C^*$-algebra $\A\rtimes X$ containing both $\A$ and a unitary representation of $X$, with a prescribed commutation rule between elements of these two sets. When a $2$-cocycle of the group (with values in the unitary group of the algebra) is also given and when ``unitary representation'' is replaced by ``projective representation'' in some suitable generalized sense, then one gets a twisted crossed product. We shall be pragmatic and introduce only the object of strict interest for our situation in a somewhat ad hoc manner. In \cite{MP1} and especially in \cite{MPR1} we give a more detailed description. The abstract theory of twisted crossed products was developed in \cite{BS}, \cite{PR1} and \cite{PR2}.

So, let us start by remarking that $X=\bR^N$ is indeed a locally compact second countable group. We shall call {\it admissible} any separable $C^*$-algebra $\A$ composed of bounded, uniformly continuous complex functions on $X$ which contains $C_0(X)$ and is invariant under translations: $a\in\A$, $x\in X$ imply $a(\cdot+x)\in\A$. Thus, for any $\h\ne 0$, one can define the continuous action of $X$ by automorphisms of $\A$: 
\begin{equation*}
\theta^\h:X\rightarrow \text{Aut}(\A), \ \ \ \left[\theta^\h_{x}(a)\right](y):=a(y+\h x).
\end{equation*} 
$\theta^\h$ is a group morphism and the maps $\ X\ni x\mapsto\theta^\h_{x}(a)\in\A$, $\ a\in\A\ $ are all continuous. Let us recall the function 
\begin{equation*}
(q,x,y)\mapsto\omega^\h_B(q;x,y):=e^{-(i/\h)\Gamma_B(<q,q+\h x,q+\h x+\h y>)},
\end{equation*}
which governs the multiplication property of the magnetic translations. It can be interpreted as a map
\begin{equation*}
\omega^\h_B:X\times X\rightarrow C(X;\bT), \ \ \ \left[\omega^\h_B(x,y)\right](q):=\omega^\h_B(q;x,y)
\end{equation*} 
with values in the set of continuous functions on $X$ taking values in the $1$-torus $\bT:=\{z\in\bC\mid\vert z\vert=1\}$. It is easy to see that $\omega^\h_B$ satisfies {\it the $2$-cocycle condition}
\begin{equation*}
\omega^\h_B(x,y)\omega^\h_B(x+y,z)=\theta^\h_x\left[\omega^\h_B(y,z)\right]\omega^\h_B(x,y+z), \ \ \ \forall x,y,z\in X,
\end{equation*}
easy to check with Stokes' Theorem, since $dB=0$. It is also {\it normalized}, i.e. 
\begin{equation*}
\omega^\h_B(x,0)=1=\omega^\h_B(0,x),\ \  \forall x\in X.
\end{equation*} 

We have shown in \cite{MPR1} how to impose conditions on $B$ in order to have a good connection between $\omega^\h_B$ and the admissible $C^*$-algebra $\A$. Let us denote by $S_\A$ the Gelfand spectrum of $\A$ (the space of characters with the pointwise convergence topology). Our assumptions on $\A$ imply that $X$ can be identified with a dense subset of the locally compact, second countable space $S_\A$. We say that a continuous function on $X$ is {\it of class} $\A$ if it extends to a continuous function on $S_\A$. The $C^*$-algebra $\A$ is unital iff $S_\A$ is compact (thus a compactification of $X$) and in this case ``continuous'' means also ``bounded''; in the non-unital case many unbounded functions are allowed. If the components $B_{jk}$ of the magnetic field are of class $\A$ the mapping $X\times X\ni(x,y)\mapsto\omega^\h_B(\cdot;x,y)\in C(S_\A;\bT)$ is well-defined and continuous with respect to the topology of uniform convergence on compact subsets of $S_\A$. (Note that $C(S_\A;\bT)$ is exactly the unitary group $\U\M(\A)$ of the multiplier algebra of $\A$.) These are the needed conditions to call $\left(\theta^\h,\omega_B^\h\right)$ {\it a twisted action of $X$ on} $\A$ and to make the quadruplet $\left(\A,\theta^\h,\omega_B^\h,X\right)$ a particular case of {\it a twisted $C^*$-dynamical system}. These are also conditions under which one can perform the construction of the twisted crossed product $C^*$-algebra that we now explain. 

Consider first the Banach space $L^1(X;\A)$ with the norm
$\parallel\vfi\parallel_1:=\int_X dx\parallel\vfi\parallel_\A$. As a rule, its
elements will be considered as functions of two variables:
$[\vfi(x)](q)\equiv\vfi(q;x)$, thus $\parallel\vfi\parallel_1=\int_X
dx\sup_{q\in X}\vert\vfi(q;x)\vert$. We can introduce an involution by
$\vfi^\diamond(q;x):=\overline{\vfi(q;-x)}$ and a composition law
\begin{equation}\label{caut}
\left(\vfi\diamond^\h\psi\right)(q;x):= 
\end{equation}
\begin{equation*}
\int_X dy\ \vfi\left(q-\frac{\h}{2}(x-y);y)\right)\psi\left(q+\frac{\h}{2}y;x-y\right)e^{-(i/\h)\Gamma_B\left(\left<q-\frac{\h}{2}x,q-\frac{\h}{2}x+\h y,q+\frac{\h}{2}x,\right>\right)}
\end{equation*} 
(we leave to the reader the task of suppressing the variable $q$ and introducing the objects $\theta^\h$ and $\omega_B^\h$ in the right places to get a more abstract version of this formula). Endowed with this structure $L^1(X;\A)$ is a Banach $^*$-algebra. 

Its envelopping $C^*$-algebra will be called {\it the twisted crossed product of $\A$ by the twisted action $\left(\theta^\h,\omega_B^\h\right)$ of} $X$. A comprehensive but awkward notation would be $\A\rtimes_{\theta^\h}^{\omega_B^\h}X$, which we abbreviate to $\gC^\h_\A$, insisting on its dependence on $\h$ and $\A$, the magnetic field $B$ being fixed. We recall that $\gC^\h_\A$ is the completion of $L^1(X;\A)$ under the $C^*$-norm 
\begin{equation*}
\pl\vfi\pl_\h:=\sup\{\pl\pi(\vfi)\pl_{B(\H)}\ \mid\pi:L^1(X;\A)\rightarrow B(\H)\ \text{representation}\}.
\end{equation*}
The main reason for $\gC^\h_\A$ to exist is the fact that its non-degenerate representations are in a one-to-one correspondence with {\it covariant representations} of the twisted $C^*$-dynamical system $\left(\A,\theta^\h,\omega_B^\h,X\right)$, i.e with triples $(\H,r,U)$, where $\H$ is a Hilbert space, $r$ is a non-degenerate representation of $\A$ and $U$ is a strongly continuous map from $X$ to the family of unitary operators on $\H$ satisfying for all $x,y\in X$ and $a\in \A$
\begin{equation}\label{covrep}
U(x)U(y)=r\left[\omega^\h_B(x,y)\right]U(x+y)\ \ \ \text{and}\ \ \ U(x)r(a)U(x)^*=r\left[\theta^\h_x(a)\right]. 
\end{equation}

We shall use this for a single case, that of {\it the Schr\"odinger covariant
  representation} $\left(L^2(X),r,U^\h_A\right)$ {\it associated to the vector
  potential} $A$ (with $dA=B$). Here $r:\A\rightarrow B[L^2(X)]$ is the usual
representation of functions in $\A$ by multiplication operators ($r(a)\equiv
a(Q)$ by a previous notation) and $U^\h_B$ has been introduced at
(\ref{uuu}). It is easy in this case to check (\ref{covrep}) and to view it as
another way to codify the commutation relations between positions and magnetic
momenta. In fact this is the root of the close connection (see below) between
$\gC^\h_\A$ and the magnetic pseudodifferential calculus sketched above. The
representation of $\gC^\h_\A$ corresponding to $\left(L^2(X),r,U^\h_A\right)$
is given (by abstract principles) by 
\begin{equation*}
\Rep^\h_A(\vfi):=\int_X dx\ r\left[\theta^\h_{x/2}(\vfi(x))\right]U^\h_A(x),
\end{equation*}
which gives for $\vfi\in L^1(X;\A) $ and $u\in L^2(X)$
\begin{equation}\label{Rep}
\left[\Rep^\h_A(\vfi)u\right](x)=\h^{-N}\int_X dy\ e^{\frac{i}{\h}\Gamma_A([x,y])}\vfi\left(\frac{x+y}{2},\frac{y-x}{\h}\right)u(y).
\end{equation}

By comparing (\ref{Rep}) with (\ref{Op}) one sees that, at least formally, $\Rep^\h_A$ and $\Op^\h_A$ are connected to each other by a partial Fourier transformation: $\Op^\h_A(f)=\Rep^\h_A\left[\bF(f)\right]$, with $\bF:=1\otimes\F$ and $(\F b)(x):=\int_{X^\star} dk\ e^{-ix\cdot k}b(k)$ whenever it makes sense. It follows that the composition laws $\circ^\h$ and $\diamond^\h$ are intertwined by $\bF$, i.e. $f\circ^\h g=\bF^{-1}\left[(\bF f)\diamond^\h (\bF g)\right]$, as can also be checked by a direct calculation. We send to \cite{MPR1} for details on the rigorous meaning of these connections in non-trivial cases. We don't need it here since actually all our verifications in Sections 5, 6 and 7 are done in the setting of twisted crossed products. One defines the $C^*$-algebra $\gB^\h_\A:=\bF^{-1}\gC^\h_\A$. On suitable dense subsets of $\gB^\h_\A$ we are entitled to use (\ref{circ}) as it stands.

\section{The main result}

In Section 2, assuming that the components of our magnetic field $B$ are
$C^\infty$ functions on $X=\bR^N$, we endowed the space $C^\infty(\Xi)_\bR$ of
real smooth functions on the phase-space $\Xi= X\times X^\star$ with a
$B$-dependent Poisson algebra structure, called $\mathcal P_B(\Xi)$.

On the other hand, in Section 3 we constructed for each $\h\in(0,1]$ a $C^*$-algebra $\gB^\h_\A$, which is the partial Fourier transform of the twisted crossed product $C^*$-algebra $\gC^\h_\A$ defined by the twisted action $\left(\theta^\h,\omega^\h_B\right)$ of $X$ on the admissible C$^*$-algebra $\A$; we had to use the assumption that the components $B_{jk}$ are functions of class $\A$.

In order to construct now a strict deformation quantization we have to study
the conditions to be imposed to the magnetic field in connection with the
choice of the Poisson subalgebra $\gA_0$ of $\mathcal P_B(\Xi)$ . 

We recall that $\gC^\h_\A$ is a $C^*$-completion of the Banach $^*$-algebra $L^1(X;\A)$; the structure depends on $\h$ and $B$. Then for any subspace $\A_0$ of $\A$ and any subspace $\S$ of $L^1(X)$, the algebraic tensor product $\A_0\odot \S$ (finite combination of simple tensors) is a subspace of $L^1(X;\A)$, thus also of $\gC^\h_\A$. The partial Fourier transformed version $\bF^{-1}\left[\A_0\odot \S\right]=\A_0\odot \F^{-1} \S$ will be a subspace of $\bF^{-1}\left[L^1(X;\A)\right]$ and, therefore, a subspace of $\gB^\h_\A$. Note that $\A_0\odot \F^{-1}\S$ is also contained in $\A\odot C_0(X^\star)$, thus it is composed of complex functions defined on the phase space $\Xi$. If one also requires that $\A_0\subset C^\infty(X)$ and $\F^{-1} \S\subset C^\infty(X^\star)$, then $\A_0\odot \F^{-1}\S\subset C^\infty(\Xi)$ and both the classical and the quantum formalisms hold on $\A_0\odot \F^{-1}\S$.
In fact several choices for $\A_0$ and $\S$ are available, their success
hanging on the assumptions we impose on the magnetic field. With severe
contraints on $B$ one hopes to quantize larger classes of classical
symbols. We shall study a simple, convenient situation; the reader could work
out other cases for himself. We define $\A^\infty:=\{a\in\A\cap
C^\infty(X)\mid\partial^\alpha a\in\A,\ \forall \alpha\in \bN^N\}$; it is a
subspace of $\A\cap BC^\infty(X)$. Take $\A_0=\A^\infty$ and $\S=\S(X)$, the Schwartz space of functions on $X$ which have rapidly decaying derivatives of any order. Then $\F^{-1}\S=\S(X^\star)$ is the Schwartz space defined on $X^\star$.

We also consider $\S(X^\star;\A^\infty)$, the space of functions $X^\star\ni p\mapsto f(p)\in\A^\infty$ such that for any $l,m\in\mathbb{N}$
$$
\|f\|_{l,m}:=\max\{\underset{p\in X^\star}{\sup}\|p^\alpha(\partial^\beta f)(p)\|_{\A}\;\mid\; |\alpha|\leq l, |\beta|\leq m\} < \infty.
$$
We remark that 
$$
\S(X^\star;\A^\infty)\subset C^\infty(\Xi)\cap\bF^{-1}\{L^1(X;\A)\}.
$$
Then we have the following evident  

\begin{Proposition}
Suppose that the components of the magnetic field $B$ belong to $\A^\infty$. Then $\S(X^\star;\A^\infty)_\bR$ is a Poisson subalgebra of $\mathcal{P}_B(\Xi)$ and a dense subset of the self-adjoint part of the abelian $C^*$-algebra $\A\otimes C_0(X^\star)$.
\end{Proposition}

We can now state

\begin{Theorem}\label{MR} {\bf (Main result)}  
Assume that the components of the magnetic field $B$ belong to $\A^\infty$. Then 
the family of injections $\left(\S(X^\star;\A^\infty)_\bR\hookrightarrow\gB^\h_\A\right)_{\h\in[0,1]}$ is a strict deformation quantization (cf. Definitions \ref{strict} and \ref{strictdefQ}).
\end{Theorem}

As seen in Section 3, one may say that $\gB^\h_{\A}$ is a $C^*$-algebra of (magnetic) pseudodifferential symbols and its represented versions $\Op^\h_A\left(\gB^\h_{\A}\right)\subset B(L^2(X))$ are $C^*$-algebras of magnetic pseudodifferential operators. It will be more convenient to work in the other realization, that of twisted crossed products. There are two reasons:\\ 1. There exist results of \cite{PR2} and \cite{Ni} on continuous fields of twisted crossed products which lead almost immediately to Rieffel's condition.\\2. In the twisted crossed product formalism one disposes of the simple norm $\pl\cdot\pl_1$, which will be very convenient in checking the axioms of von Neumann and Dirac.

Thus we state now a variant of Theorem \ref{MR}; these two results are
equivalent by the isomorphisms defined by the partial Fourier
transformation. We need first to rewrite the magnetic Poisson structure. On $\S(X;\A)$ (obvious definition) we set by transport of structure
\begin{equation*}
\vfi\diamond^0\psi:=\bF\left[(\bF^{-1}\vfi)(\bF^{-1}\psi)\right]\ \ \text{and}\ \ \{\vfi,\psi\}^B:=\bF\left[\{\bF^{-1}\vfi,\bF^{-1}\psi\}_B\right],\ \ \ \vfi,\psi\in\S(X;\A).
\end{equation*}
A simple direct calculation gives
\begin{equation}\label{unu}
(\vfi\diamond^0\psi)(q;x)=\int_X dy\ \vfi(q;y)\psi(q;x-y);
\end{equation}
$\diamond^0$ is poinwise multiplication in the first variable and convolution
in the second. Slightly more effort is needed to prove that 
\begin{equation*}
\{\vfi,\psi\}^B=-i\sum_{j=1}^N\left[(Q^{(2)}_j\vfi)\diamond^0(\pt^{(1)}_j\psi)-(\pt^{(1)}_j\vfi)\diamond^0(Q^{(2)}_j\psi)\right]-
\end{equation*}
\begin{equation}\label{doi}
-\sum_{j,k=1}^N B_{jk}(\cdot)(Q^{(2)}_j\vfi)\diamond^0(Q^{(2)}_k\psi),
\end{equation}
where $\left(Q^{(2)}_j\rho\right)(q;x):=x_j\rho(q;x)$ and $(\pt^{(1)}_j\rho)(q;x)=\frac{\pt}{\pt q_j}\rho(q;x)$.

Let us denote by $C^*(X)$ the group $C^*$-algebra of $X$; it is the envelopping $C^*$-algebra of $L^1(X)$, the convolution Banach $^*$-algebra of $X$. It is isomorphic to $C_0(X^\star)$ by an extension of the Fourier transformation; thus the spectrum of $C^*(X)$ is homeomorphic to $X^\star$. Note that the twisted crossed product $\gC^\h_\A=\A\rtimes_{\theta^\h}^{\omega^\h_B} X$ collapses to $\A\otimes C^*(X)$ for $\h=0$.
\begin{Proposition}
Suppose that the components of the magnetic field $B$ belong to $\A^\infty$; then the vector space $\S(X;\A^\infty)_\bR=\bF\left[\S(X^\star;\A^\infty)\right]_{\mathbb{R}}$ is a Poisson algebra for the composition laws (\ref{unu}) and (\ref{doi}). It is also dense in the self-adjoint part of the abelian $C^*$-algebra $\A\otimes C^*(X)$.
\end{Proposition}

The partial Fourier transformed version of our Main Result reads

\begin{Theorem} Assume that the components of the magnetic field $B$ belong to
  $\A^\infty$; then the family of injections
  $\left(\S(X;\A^\infty)_\bR\hookrightarrow(\gC^\h_\A)_\bR\right)_{\h\in[0,1]}$ is a strict deformation quantization (cf. Definitions \ref{strict} and \ref{strictdefQ}).
\end{Theorem}
The completeness condition is obvious: $\S(X;\A^\infty)$ is dense in
$\left(L^1(X;\A),\pl\cdot\pl_1\right)$, $L^1(X;\A)$ is dense in
$\left(\gC^\h_\A,\pl\cdot\pl_\h\right)$ and one has $\pl\cdot\pl_1\ \le\
\pl\cdot\pl_\h$. The conditions of Definition \ref{strictdefQ} are also
clearly satisfied. We still have to verify the conditions (a), (b) and
(c) of Definition \ref{strict}. This will be done in the next sections. 

{\bf Remark.} It would be in the spirit of many works in strict deformation quantization to consider only the case $\A=C_0(X)$. Since in this case $\gC^\h_{C_0(X)}$ is isomorphic to $K\left[L^2(X)\right]$, the $C^*$-algebra of all compact operators on $L^2(X)$ (cf. \cite{MPR1}, Proposition 2.17 (b)), in fact one works with a field of $C^*$-algebras with two types of fibers: $C_0(\Xi)$ for $\h=0$ and $K\left[L^2(X)\right]$ for $\h\ne 0$. We think that both the twisted crossed product $\gC^\h_\A$ and the pseudodifferential formalism are useful for arbitrary, admissible $\A$. In \cite{MPR2} it is shown how to calculate the essential spectrum and how to get localization results for generalized Schr\"odinger operators with anisotropic potentials and magnetic fields. The anisotropy is taken into account by the abelian algebra $\A$ and exploiting the structure of its spectrum is the key of the proofs.

{\bf Remark.} Let us point out that if the spectrum of $\A$ is compact (and
that is always the case in the applications to quantum Hamiltonians, where we
expect $\A$ to have a unit), then the components of the magnetic field $B$
being of type $\A$ evidently imply that they are bounded and uniformly
continuous. Thus, in this case the requirement that the components of $B$ are
of class $\A^\infty$ (i.e. they are of class $C^\infty(X)$ and together with
all their derivatives admit continuous extensions to the spectrum of $\A$) is
rather optimal. If we allow the spectrum of $\A$ to be noncompact, then we can
allow unbounded magnetic fields with components of class $\A$ but we have to
replace $\A^\infty$ with $\A^\infty_{\rm c}$ the subalgebra of elements of $\A^\infty$ that have compact support (with respect to the spectrum of $\A$).

An important technical ingredient in our proof relies on a result saying roughly that, under certain conditions, the twisted crossed product of a group with the sectional algebra of a $C^*$-bundle is the sectional algebra of a $C^*$-bundle of twisted crossed products. This can be found in \cite{PR2} and \cite{Ni}; techniques of \cite{Ri1} and \cite{Bl} are also relevant here. For us the most convenient reference is \cite{Ni}, from which we quote slightly reformulated the definition and the result below. 
\begin{Definition}
{\rm A continuous $C^*$-bundle} is a triple $\mathbf A=\left(I,\{\A^\h\}_{\h\in I},\Gamma_0(\mathbf A)\right)$, where $I$ is a Hausdorff, locally compact space, $\A^\h$ is a $C^*$-algebra with norm $\pl\cdot\pl_\h$ and $\Gamma_0(\mathbf A)$ a $C^*$-algebra of sections such that:

(i) For any $\h\in I$, $\ \{F(\h)\mid F\in\Gamma_0(\mathbf A)\}=\A^\h$.

(ii) For any $F\in\Gamma_0(\mathbf A)$, the map $\ \h\mapsto\pl F_\h\pl_\h$ belongs to $C_0(I)$. 

(iii) $\Gamma_0(\mathbf A)$ is a $C_0(I)$-module: if $F\in\Gamma_0(\mathbf A)$ and  $\nu\in C_0(I)$, then $\nu F$ (defined pointwise) also belongs to $\Gamma_0(\mathbf A)$.
\end{Definition}
In fact the arguments in \cite{PR2} show that the separability condition in their definition of the twisted crossed-product is needed only in studying the structure of the group of cocycles. Thus, for our developments of the functional calculus with magnetic fields, we can consider a slightly general definition for twisted croosed-products by eliminating the separability condition, and as the proof in \cite{Ni} is still valid, we have in fact the theorem cited below.
\begin{Theorem}\label{Nilsen} {\rm [Nielsen 1996 \cite{Ni}]}
Let $\mathbf A$ be a continuous $C^*$-bundle such that $\Gamma_0(\mathbf A)$ is separable. Let $(\Theta,\Omega)$ be a twisted action of an amenable, second countable locally compact group $X$ on $\Gamma_0(\mathbf A)$ by $C_0(I)$-automorphisms. Then there exists a continuous $C^*$-bundle $\mathbf C=\left(I,\{\C^\h\}_{\h\in I},\Gamma_0(\mathbf C)\right)$ such that: 

(i) For any $\h\in I$, $\ \C^\h=\A^\h\rtimes^{w^\h}_{t^h} X$, where $t_x^\h:X\rightarrow \rm{Aut}(\A^\h)$, $t_x^\h[F(\h)]:=[\Theta_x(F)](\h)$, $\forall x\in X$, $\forall F\in\Gamma_0(\mathbf A)$ and $w^\h:X\times X\rightarrow \U\M\left(\A^\h\right)$, $w^\h(x,y):=[\Omega(x,y)](\h)$, $\forall x,y\in X$.

(ii) The map $\left[(\chi\Phi)(\h)\right](x)=[\Phi(x)](\h)$, $\forall \h\in I, \forall x\in X$ extends to an isomorphism $\chi:\Gamma_0(\mathbf A)\rtimes_\Theta^\Omega X\rightarrow \Gamma_0(\mathbf C)$ such that for every $\Phi\in L^1(X;\Gamma_0(\mathbf A))$ one has $(\chi\Phi)(\h)\in L^1(X;\A^\h)$.
\end{Theorem}

\section{Rieffel's condition}

We are placed in the framework of Section 3. We start by constructing a twisted action on a large $C^*$-algebra, consisting of functions depending both on the variables $\h\in[0,1]$ and $q\in X\equiv \bR^N$. The same strategy has been used in \cite{Be} for the rotation algebras (which are also twisted crossed products) in order to explore the regularity of the spectrum of certain finite-difference operators, the parameter $\h$ being replaced there by the strength of a (discrete) magnetic field. 

We consider first the $C^*$-bundle $\mathbf A=\left(I,\{\A^\h\}_{\h\in I},\Gamma_0(\mathbf A)\right)$, where $I:=[0,1]$ is compact, $\A^\h:=\A$ (our admissible $C^*$-algebra) for all $\h$ and $\Gamma_0(\mathbf A):=C(I;\A)$. One checks easily that $\mathbf A$ is indeed a continuous $C^*$-bundle. Note that the Gelfand spectrum of the $C^*$-algebra $C(I;\A)$ is homeomorphic to $I\times S_\A$, where $S_\A$ is the spectrum of $\A$. Recalling the twisted actions $\{\left(\theta^\h,\omega^\h_B\right)\mid \h\in I\}$ of Section 3, one defines for all $\h\in I$, $q,x,y\in X$ and $F\in C(I;\A)$:
\begin{equation}\label{Theta}
\Theta:X\rightarrow\text{Aut}[C(I;\A)],\ \ \ \left(\Theta_x F\right)(\h):=\theta^\h_x[F(\h)],
\end{equation}
\begin{equation}
\Omega_B:X\times X\rightarrow C(I\times S_\A;\bT),\ \ \ \left[\Omega_B(x,y)\right](\h,q):=\omega_B^\h(q;x,y).
\end{equation}
By using notations as $[F(\h)](q)\equiv F(\h,q)$ (the elements of
$C(I;\A)\cong C(I\times S_\A)$ may be seen as functions on $I\times X$),
(\ref{Theta}) can be rewritten $\left(\Theta_x F\right)(\h,q)=F(\h,q+\h x)$.

The group $X=\bR^N$, being abelian, is amenable. Then it is easy to verify that $\left(C(I;\A),\Theta,\Omega_B,X\right)$ is a twisted $C^*$-dynamical system and that $\Theta_x(\nu F)=\nu\Theta_x(F)$ for all $x\in X$, $\nu\in C(I)$ and $F\in C(I;\A)$, as required by Theorem \ref{Nilsen}.

To apply Theorem \ref{Nilsen}, one must compute first the twisted actions
$\{\left(t^\h,w^\h\right)\mid\h\in I\}$ associated to
$\left(\Theta,\Omega_B\right)$. It easily comes out that $t^\h=\theta^\h$ and
$w^\h=\omega^\h_B$; just use the explicit formulae. Thus the $C^*$-algebras
$\C^\h$, the fibers of the continuous $C^*$-bundle $\mathbf C$, coincide
(respectively) with the $C^*$-algebras $\gC^\h_\A$ defined in Section 3. To
show that the map $\h\mapsto\pl\vfi\pl_\h$ is continuous for any
$\vfi\in\A^\infty\odot\S(X)$, by the axiom (ii) of a continuous $C^*$-bundle,
one has just to prove that any element $\varphi\in S(X;\A^\infty)$ defines a
(constant) section belonging to $\Gamma_0(\mathbf C)$. This is obvious even
for $\varphi\in L^1(X;\A)$, since the isomorphism $\chi$ just intertwins the
variables $\h$ and $x$.

\section{The von Neumann condition}

We have to show that, for fixed $\vfi,\psi\in\S(X;\A^\infty)_\bR$ we have
\begin{equation*}
\lim_{\h\rightarrow 0}\pl\frac{1}{2}\left(\vfi\diamond^\h\psi+\psi\diamond^\h\vfi\right)-\vfi\diamond^0\psi\pl_\h=0.
\end{equation*}
The operations $\diamond^\h$ and $\diamond^0$ are defined, respectively, at (\ref{caut}) and (\ref{unu}).
Taking into account that $\pl\cdot\pl_\h\ \le\ \pl\cdot\pl_1$ and by the triangle inequality, it is enough to prove 
\begin{equation}\label{cvN}
\lim_{\h\rightarrow 0}\pl\vfi\diamond^\h\psi-\vfi\diamond^0\psi\pl_1=0.
\end{equation}

By standard arguments one can approach any function in $L^1(X;\Gamma_0(\mathbf{A})$, in $L^1$-norm, with a continuous function with compact support, and show that (using the notations of Theorem \ref{Nilsen})
$$
\chi[L^1(X;\Gamma_0(\mathbf{A}))]= C(I;L^1(X,\A)).
$$

For any function $\varphi\in L^1(X;\A)$ let us denote by $\varphi_\circ\in C(I;L^1(X;\A)$ the constant function $\varphi_\circ(\h):=\varphi,\;\forall\h\in I$ and by $\tilde{\varphi}\in L^1(X;C(I;\A))$ the function taking constant values $[\tilde{\varphi}(x)](\h):=\varphi(x),\;\forall\h\in I$. Let us also remark that any constant function in $C(I;L^1(X;\A)$ is of the form $\varphi_\circ=\chi[\tilde{\varphi}]$ for some $\varphi\in L^1(X;\A)$. We denote by $\sharp$ the product in $\Gamma_0(\mathbf{A})\rtimes^\Omega_\Theta X$ and by $\diamond$ the product in $\Gamma_0(\mathbf{C})$. Then for fixed $\varphi,\psi\in\S(X;\A^\infty)_\bR\subset L^1(X;\A)$, one has
$$
[\varphi_\circ\diamond\psi_\circ](\h)=
\chi[\tilde{\varphi}\;\sharp\;\tilde{\psi}](\h)
$$
and thus $\varphi_\circ\diamond\psi_\circ\in C(I;L^1(X;\A))$. As a direct consequence, if we set $\boldsymbol{\Psi}(\varphi,\psi):=\varphi_\circ\diamond\psi_\circ-
(\varphi\diamond^0\psi)_\circ\in C(I;L^1(X;\A)$ we get
$$
\underset{\h\rightarrow 0}{\lim}\|\boldsymbol{\Psi}(\varphi,\psi)(\h)-
\boldsymbol{\Psi}(\varphi,\psi)(0)\|_1=0.
$$
But 
$$
\boldsymbol{\Psi}(\varphi,\psi)(0)=(\varphi_\circ\diamond\psi_\circ)(0)-
(\varphi\diamond^0\psi)_\circ(0)=\varphi\diamond^0\psi-\varphi\diamond^0\psi=0
$$
so that
$$
\|\varphi\diamond^\h\psi-\varphi\diamond^0\psi\|_\h=
\|\boldsymbol{\Psi}(\varphi,\psi)(\h)\|_\h\leq
\|\boldsymbol{\Psi}(\varphi,\psi)(\h)\|_1\underset{\h\rightarrow 0}{\rightarrow}0
$$
and thus we have got the von Neumann condition.

An alternative direct proof by elementary arguments can be given by studying the exponent of the cocycle appearing in the definition of $\diamond^\h$:
$$
-\frac{i}{\h}\Gamma_B(<q-\frac{\h}{2}x,q-\frac{\h}{2}x+\h y,q+\frac{\h}{2}x>)=-\frac{i}{\h}\int\limits_{<q-\frac{\h}{2}x,q-\frac{\h}{2}x+\h y,q+\frac{\h}{2}x>}B.
$$
Consider a parametrization for the triangle $<q-\frac{\h}{2}x,q-\frac{\h}{2}x+\h y,q+\frac{\h}{2}x>$:
$$
<q-\frac{\h}{2}x,q-\frac{\h}{2}x+\h y,q+\frac{\h}{2}x>\;=\;\kappa_{q,(x,y)}[\Delta_2]
$$
where $\Delta_2:=\{(t,s)\in\mathbb{R}^2\;\mid\;0\leq t \leq 1,\; 0\leq s\leq t\}$ and
$$
\kappa_{q,(x,y)}(t,s):=q-\frac{\h}{2}x+t\h y+s\h (x-y)=q+(s-\frac{1}{2})\h x+(t-s)\h y.
$$
Thus, by denoting $e_t$ and $e_s$ the tangent vectors corresponding to the two coordinate functions of $\mathbb{R}^2$, we have
$$
\int_{\kappa[\Delta_2]}B=\int_{\Delta_2}\kappa^*B=\int_0^1dt\int_0^tds\left(\kappa^*B\right)(e_t,e_s).
$$
An obvious calculation gives
$$
(\kappa^*B)(e_t,e_s)=\sum_{j,k}B\left(\kappa_{q,(x,y)}(t,s)\right)\frac{\partial\kappa_j}
{\partial t}\frac{\partial\kappa_k}{\partial s}
$$
and thus we have
\begin{equation*}
-\frac{i}{\h}\Gamma_B(\kappa_{q,(x,y)}[\Delta_2])=-i\h\sum_{j,k}y_j(x-y)_k\int_0^1dt\int_0^tds
B_{jk}\left[q+(s-\frac{1}{2})\h x+(t-s)\h y\right]
\end{equation*}
\begin{equation}\label{xi}
=:-i\h \Omega_B(q,x,y;\h).
\end{equation}
Now let us come back to (\ref{cvN}) and estimate
$$
\pl\vfi\diamond^\h\psi-\vfi\diamond^0\psi\pl_1\leq
$$
\begin{equation*}
\leq\int_Xdx\int_Xdy\; \underset{q\in X}{\sup}\Big|
\vfi\left(q-\frac{\h}{2}(x-y);y\right)\psi\left(q+\frac{\h}{2}
y;x-y\right)e^{-i\h\Omega_B(q,x,y;\h)}\;-
\end{equation*}
\begin{equation}\label{I1}
-\vfi(q;y)\psi(q;x-y)\Big|.
\end{equation}
It is easy to see that the integral is bounded by the expression
$$
2\left(\underset{q\in X}{\sup}\;\underset{y\in X}{\sup}<y>^m\vfi(q;y)\right)\left(\underset{q\in X}{\sup}\;\underset{x\in X}{\sup}<x>^m\psi(q;x)\right)
\left(\int_Xdx<x>^{-m}\right)^2,
$$
that is finite and not depending of $\h$ for any $m>N$.
On the other hand the integrand in (\ref{I1}) is convergent pointwise to zero
when $\h$ goes to $0$, as can be seen after writing the inequality
$$
\underset{q\in X}{\sup}\left|
\vfi\left(q-\frac{\h}{2}(x-y);y\right)\psi\left(q+\frac{\h}{2} y;x-y\right)e^{-i\h\Omega_B(q,x,y;\h)}-\vfi(q;y)\psi(q;x-y)\right|\leq
$$
$$
\leq\left(\underset{q\in X}{\sup}\left|\vfi\left(q-\frac{\h}{2}(x-y);y\right)-
\vfi(q;y)\right|\right)\left(\underset{q\in X}{\sup}\left|\psi\left(q+\frac{\h}{2}y;x-y\right)\right|\right)
+
$$
$$
+\left(\underset{q\in X}{\sup}\left|\vfi(q;y)\right|\right)\left(\underset{q\in X}{\sup}\left|\psi\left(q+\frac{\h}{2}y;x-y\right)-
\psi(q;x-y)\right|\right)+
$$
$$
+\left(\underset{q\in X}{\sup}\left|\vfi(q;y)\psi(q;x-y)\right|\right)
\left(\underset{q\in X}{\sup}\left|e^{-i\h\Omega_B(q,x,y;\h)}-1\right|\right).
$$
For the first two lines we use the fact that $\vfi$ and $\psi$ belong to
$\S(X;\A^\infty)\subset BC^\infty(X\times X)$. For the third one, the hypothesis that the components of the magnetic field are in $\A^\infty\subset BC^\infty(X)$ implies that
for any $(x,y)\in X^2$ we have $\underset{q\in
  X}{\sup}|\Omega_B(q,x,y;\h)|\leq C(x,y)$ uniformly in $\h$. Thus $\;\underset{q\in X}{\sup}|\exp\{-i\h \Omega(q,x,y;\h)\}-1|\;$ converges to $0$ for
$\h\rightarrow 0$.

\section{The Dirac condition}

We need only to prove that the following convergence holds:
\begin{equation}\label{cD}
\left\|\frac{1}{i\h}\left(\vfi\diamond^\h\psi-\psi\diamond^\h\vfi\right)
-\{\vfi,\psi\}_B\right\|_1\underset{\h\rightarrow 0}{\longrightarrow} 0.
\end{equation}
For that we shall need the exact form of the first order term in $\h$ of $\;\vfi\diamond^\h\psi-\psi\diamond^\h\vfi$. We use Taylor developments
\begin{equation*}
\vfi\left(q-\frac{\h}{2}(x-y);y\right)=\vfi(q;y)-\frac{\h}{2}\sum_{j=1}^N (x_j-y_j)\int_0^1 ds\ \left(\pt^{(1)}_j\vfi\right)\left(q-s\frac{\h}{2}(x-y);y\right)
\end{equation*}
and
\begin{equation*}
\psi\left(q+\frac{\h}{2}y;x-y\right)=\psi(q;x-y)+\frac{\h}{2}\sum_{j=1}^N y_j\int_0^1 ds\ \left(\pt^{(1)}_j\psi\right)\left(q+s\frac{\h}{2}y;x-y\right).
\end{equation*}
For $\rho\in\S(X;\A^\infty)$, $z\in X$ and $(q,x)\in X^2$ we shall use the
notation ${\nabla}^{(1)}\rho$ for the gradient with respect to the
first variable in $X\times X$ and set
\begin{equation}\label{defL}
(\mathfrak{L}_z^{\pm}\rho)(q;x):=\frac{1}{2}z\cdot\int_0^1ds
\;\left({\nabla}^{(1)}\rho\right)\left(q\pm s\frac{\h}{2}z;x\right).
\end{equation}
Note that we have
\begin{equation}\label{ultima}
\underset{\h\rightarrow 0}{\lim}(\mathfrak{L}_z^{\pm}\vfi)(q;x)
=\frac{1}{2}z\cdot\left({\nabla}^{(1)}\vfi\right)(q;x).
\end{equation}
Then the Taylor developments above read
$$
\vfi\left(q-\frac{\h}{2}(x-y);y\right)=
\vfi(q;y)-\h(\mathfrak{L}_{x-y}^{-}\vfi)(q;y)
$$
and
$$
\psi\left(q+\frac{\h}{2}y;x-y\right)=\psi(q;x-y)+
\h(\mathfrak{L}_y^{+}\vfi)(q;x-y).
$$

Moreover, the assumption that $B_{jk}\in\A^\infty\subset BC^\infty(X)$ implies that for any $x,y\in X$, the function $X\times [0,1]\ni(q,\h)\mapsto\Omega_B(q,x,y;\h)\in\mathbb{R}$ is bounded and uniformly continuous. Thus, if we denote
\begin{equation}\label{restD}
R_B(q,x,y;\h):=\frac{1}{\h}\left(e^{-i\h\Omega_B(q,x,y;\h)}-e^{-i\h\Omega_B(q,x,y,0)}\right),
\end{equation}
we get $\ \underset{\h\rightarrow 0}{\lim}|R_B(q,x,y;\h)|=0$.

Putting everything together we obtain
$$
\left(\vfi\diamond^\h\psi-\psi\diamond^\h\vfi\right)(q;x)=
$$
$$
=\int_Xdy\;\left[(1-\h\mathfrak{L}^-_{x-y})\vfi\right](q;y)\left[(1+
\h\mathfrak{L}^+_y)\psi\right](q;x-y)\left[e^{-i\h\Omega_B(q,x,y;0)}+\h R_B(q,x,y;\h)\right]-
$$
$$
-\int_Xdy\;\left[(1-\h\mathfrak{L}^-_{x-y})\psi\right](q;y)\left[(1+
\h\mathfrak{L}^+_y)\vfi\right](q;x-y)\left[e^{-i\h\Omega_B(q,x,y;0)}+\h R_B(q,x,y;\h)\right]=
$$
$$
=\int_Xdy\;\vfi(q,y)\psi(q;x-y)\left[e^{-i\h\Omega_B(q,x,y;0)}-
e^{-i\h\Omega_B(q,x,x-y;0)}\right]+
$$
$$
+\h\int_Xdy\;\vfi(q;y)\left[(\mathfrak{L}^+_y+\mathfrak{L}^-_y)
\psi\right](q;x-y)-
\h\int_Xdy\;\left[(\mathfrak{L}^+_{x-y}+\mathfrak{L}^-_{x-y})
\vfi\right](q;y)\psi(q;x-y)+o(\h)
$$
where, for obtaining the second identity, we have changed an integration
variable from $y$ to $x-y$. By using (\ref{ultima}) and some simple arguments
we get
$$
\left(\vfi\diamond^\h\psi-\psi\diamond^\h\vfi\right)(q;x)=
$$
$$
=-i\h\sum_{j,k}B_{jk}(q)\int_Xdy\;y_j\vfi(q;y)\;(x-y)_k\psi(q;x-y)+
$$
$$
+\h\int_Xdy\;\left[y\vfi(q;y)\cdot({\nabla}^{(1)}\psi)(q;x-y)-
({\nabla}^{(1)}\vfi)(q;y)\cdot(x-y)\psi(q;x-y)\right]+o(\h),
$$
 The result is now straightforward by the explicit
form of the bracket $\{\cdot,\cdot\}^B$ and of the composition law $\diamond^0$.

\enddocument